\documentclass[12pt]{amsart}
\usepackage{amsthm}
\usepackage{amsmath}
\usepackage[utf8]{inputenc}
\usepackage[top=72pt,bottom=96pt,left=72pt,right=72pt]{geometry}
\def\ad{\mathrm{ad}}

\def\Aut{\mathrm{Aut}}
\def\C{\mathbb{C}}
\def\e{\varepsilon}
\def\End{\mathrm{End}}
\def\g{\mathfrak{g}}
\def\Gr{\mathrm{Gr}}
\def\h{\mathfrak{h}}
\def\H{\mathbb{H}}

\def\hyper{\mathrm{hyper}}
\def\id{\mathrm{id}}
\def\Im{\mathrm{Im}}

\def\k{\mathfrak{k}}

\def\m{\mathfrak{m}}
\def\P{\mathrm{P}}
\def\pfill{\vskip3pt plus2pt minus2pt\noindent}
\def\pr{\mathrm{pr}}
\def\R{\mathbb{R}}
\def\Re{\mathrm{Re}}
\def\Ric{\mathrm{Ric}}
\def\S{\mathrm{Sym}}

\def\SO{\mathrm{SO}}
\def\sp{\mathfrak{sp}}
\def\Sp{\mathrm{Sp}}
\def\span{\mathrm{span}}
\def\Spin{\mathrm{Spin}}
\def\SU{\mathrm{SU}}
\def\t{\mathfrak{t}}

\def\<#1,#2>{\langle\,#1,\,#2\,\rangle}
\newtheorem{Lemma}{Lemma}[section]
\newtheorem{Theorem}[Lemma]{Theorem}
\newtheorem{Corollary}[Lemma]{Corollary}
\newtheorem{Proposition}[Lemma]{Proposition}
\theoremstyle{definition}

\newtheorem{Definition}[Lemma]{Definition}
\title{On Quaternionic Bisectional Curvature}
\author{Oscar Macia, Uwe Semmelmann \& Gregor Weingart}
\address{Oscar Macia\\
 Departamento de Matemáticas\\
 Facultad de Ciencias Matemáticas\\
 Universidad de Valencia E.~G.\\
 Av.~Vicente Andrés Estellés 1\\
 46100 Burjassot, Valencia, Spain}
\email{oscar.macia@uv.es}
\address{Uwe Semmelmann\\
 Institut für Geometrie und Topologie\\
 Fachbereich Mathematik\\
 Universität Stuttgart\\
 Pfaffenwaldring 57\\
 70569 Stuttgart, Germany}
\email{uwe.semmelmann@mathematik.uni-stuttgart.de}
\address{Gregor Weingart\\
 Instituto de Matemáticas\\
 Universidad Nacional Autónoma de México\\
 Avenida Universidad s/n\\
 Colonia Lomas de Chamilpa\\
 62210 Cuernavaca, Morelos, Mexico}
\email{gw@matcuer.unam.mx}
\date{\today}
\begin{document}
\begin{abstract}
 In this article we study the concept of quaternionic bisectional
 curvature introduced by B.~Chow and D.~Yang in \cite{chow} for
 quaternion--Kähler manifolds. We show that non--negative quaternionic
 bisectional curvature is only realized for the quaternionic projective
 space $\H\P^n$. We also show that all symmetric quaternion--Kähler
 manifolds different from $\H\P^n$ admit quaternionic lines of negative
 quaternionic bisectional curvature. In particular this implies that
 non--negative sectional curvature does not imply non--negative quaternionic
 bisectional curvature. Moreover we give a new and rather short proof of
 a classification result by A.~Gray on compact Kähler manifolds of
 non--negative sectional curvature.
\end{abstract}
\maketitle
\section{Introduction}
 A Riemannian manifold $(\,M^{4n},\,g\,)$ is said to be quaternion--Kähler,
 if its Riemannian holonomy group is contained in $\Sp(n)\cdot\Sp(1)$.
 Equivalently, if  locally there exists a family of almost complex structures
 $\{I,J,K\},$ compatible with the Riemannian metric $g$, satisfying the
 quaternion relations $I^2=J^2=K^2=IJK=-\id$ and such that 
 $I, J, K$ span a  globally defined parallel subbundle of
 $\End(TM)$. 

 Quaternion--Kähler manifolds  enjoy remarkable geometric properties. In
 particular they are Einstein and they all posses a twistor space fibration
 $ZM \longrightarrow M$ with $\C\P^1$ fibres. 
 Actually $ZM$ is defined as the unit sphere bundle of the rank--$3$ vector bundle  
 spanned by the locally defined almost complex structures. 
 In the case
 of positive scalar curvature, i.e.~for so--called  {\it positive}
 quaternion--Kähler manifolds, it is known that the twistor space $ZM$
 is Kähler--Einstein, Fano, and comes equipped with a complex contact
 structure  (\cite{lebrun1}, \cite{salamon}). The only positive
 quaternion--Kähler manifolds known to date are the Wolf spaces \cite{wolf},
 including $\H\P^n,\;\Gr_2(\C^{n+4}),$ and $G_2/SO(4)$ as prominent
 examples. They are all symmetric spaces, grouped in three infinite families
 and five exceptional cases. The typical example of the twistor fibration is
 $\C\P^{2n+1}\longrightarrow\H\P^n$.

 The central long--standing conjecture in the classification problem of
 compact quaternion--Kähler manifolds is credited to  LeBrun \& Salamon
 \cite{lebrun}, and states that  {\it any positive quat\-ernion--Kähler
 manifold is symmetric, thus one of the Wolf spaces}. Despite many advances
 and partial results, this conjecture has been proved completely only for
 $1\,\leq\,n\,\leq\,4$ \cite{hitchin, poon, buczynski}.

 Inspired by results in Kähler geometry, by imposing an additional hypothesis
 on the  sectional curvature, a weaker conjecture may be obtained:
 {\it any positive quaternion--Kähler manifold of non--negative sectional
 curvature is a symmetric space}. The analogous conjecture is a known
 theorem in the realm of Kähler manifolds. Indeed, A.~Gray proved in
 \cite{gray} the following theorem:  {\it Let $M$ be a compact Kähler
 manifold with non--negative sectional curvature and constant scalar
 curvature, then $M$ is locally symmetric}. If $M$ is in addition simply 
 connected it has to be a Hermitian symmetric space. 
 Gray's proof is presented and slightly improved in \cite{HM} with
 an application to an eigenvalue estimate of the Laplace operator
 on the unit sphere bundle and similarly in \cite{SS} with an extension
 to nearly Kähler manifolds. We will begin our article with a new and 
 rather short proof of Gray's theorem. Our proof is based on the ideas of
 Gray but takes a rather different approach.

 For quaternion--Kähler manifolds even the weaker conjecture is still
 completely open. We only have the result of  Berger \cite{berger} that a
 quaternion--Kähler manifold of positive sectional curvature has to be
 isometric to $\H\P^n$. In the same
 article Berger proved that a Kähler manifold of positive sectional 
 curvature has to be isometric to $\C\P^n$

 For Kähler manifolds the concept of sectional curvature is further
 refined to {\it holomorphic bisectional curvature} and a fundamental
 theorem of N.~Mok and J.-Q.~Zhong \cite{Mok} states that any compact
 Kähler--Einstein manifold of non--negative holomorphic bisectional
 curvature and positive Ricci curvature has to be isometric to a
 compact Hermitian symmetric space.

 In the same spirit, a quaternion--Kähler analogue, {\it quaternionic
 bisectional curvature}, has been proposed by B.~Chow and D.~Yang in
 \cite{chow}. However the parallel between Kähler and quaternion--Kähler
 manifolds fails to realize as non--negative quaternionic bisectional
 curvature does not characterize quaternion--Kähler symmetric spaces.
 We will explain in our article, that the only quaternion--Kähler manifold
 satisfying this condition is $\H\P^n$ (see Corollary \ref{nnq}).
 Explicitly  we will show that on all Wolf spaces different from $\H\P^n$,
 there always exist quaternionic lines producing negative quaternionic
 bisectional curvature (see Theorem \ref{main}). This in particular shows,
 contrary to a remark in \cite{chow}, that non--negative sectional curvature
 does not imply non--negative quaternionic bisectional curvature. In fact
 \cite{chow} is written in a rather misleading way. The article makes the
 impression that it has proved that quaternion--Kähler manifolds of
 non--negative sectional curvature have to be symmetric. In this way
 the conclusion was cited e.~g.~by Brendle and Schoen in \cite{BS} for
 the proof of Proposition 11. For the first time it was remarked in
 \cite{MN} that there is a problem with \cite{chow} and that the
 classification problem for quaternion--Kähler manifolds on
 non--negative sectional curvature is still widely open.
\section{Kähler manifolds on non--negative curvature}
 The aim of this section is to give a new and very short proof of the result by
 A.~Gray mentioned in the introduction. We start with recalling a few facts
 from  \cite{GU1} on the curvature operator $q(R)$ and the standard
 Laplacian $\Delta_\rho$. Two well-known important results on these
 two objects will be the basis for our proof of the Gray result.
\subsection{The curvature endomorphism}
 Let $(M^m, g)$ be a Riemannian manifold and $EM$ some geometric vector
 bundle over $M$, i.e.~a vector bundle associated to the frame bundle by
 some representation $\rho: O(m) \rightarrow \Aut (E)$. We denote with
 $\nabla$ the connection on sections of $EM$ induced by the Levi--Civita
 connection of $g$. Its  curvature is defined as $R_{X, Y} = \nabla^2_{X, Y}
 - \nabla^2_{Y, X}$ for any tangent vectors $X, Y$. In this situation it
 is convenient to introduce the curvature endomorphism $q(R) \in \End (EM)$
 defined by
 $$
  q(R)
  \;\;=\;\;
  \frac12\,\sum_{\mu,\,\nu} (E_\mu\wedge E_\nu)_* \circ R_{E_\mu,\,E_\nu}
 $$
 where $\{\,E_\mu\,\}$ is an orthonormal basis and the star denotes the
 differential of the representation $\rho$ defining $E$ (see \cite{GU1},
 Section 2.4, for further details). The introduction of the curvature
 endomorphism $q(R)$ is motivated by the definition of the standard Laplace
 operator $\Delta_\rho = \nabla^*\nabla + q(R)$, which acts as a Laplace
 type operator on sections of $EM$. The operator $\Delta_\rho$ has several
 nice properties. It coincides with the Hodge Laplacian $\Delta\,=\,dd^*
 \,+\,d^*d$ on $k$--forms, i.e.~for $E\,=\,\Lambda^kT$ and similarly with
 the Lichnerowicz Laplacian on symmetric tensors, i.e.~for $E\,=\,\S^kT$.
 Moreover, if $EM$ is a homogeneous vector bundle over a symmetric space
 $M\,=\,G/K$ the action of the standard Laplace operator $\Delta_\rho$ on
 sections of $EM$ coincides with the action of the Casimir operator of $G$
 with respect to the left regular representation. For this article perhaps
 the most important property of $\Delta_\rho$ is that it commutes with
 parallel bundle maps. Note that the same is true for the curvature
 endomorphism $q(R)$.

 The curvature endomorphism is a symmetric endomorphism, linear in the
 curvature $R$, i.e.~it is pointwise diagonalisable. We will write $q(R)
 \,\geq\,0$, if all its eigenvalues are non--negative. It turns out that
 for a general tensor bundle non--negative sectional curvature does not
 imply that $q(R)$ is non--negative. As a remarkable fact we have however that
 this is true on symmetric tensors (see e.g. \cite{AU1}, proof of Proposition
 6.6).

 \begin{Proposition}\label{curv11}
  Let $(\,M,\,g\,)$ be a Riemannian manifold with non--negative sectional
  curvature then $q(R)\,\geq\,0$ as a symmetric endomorphism of $\S^kTM$.
 \end{Proposition}

 \proof
 The details of the proof can be found in \cite{AU1}. However, for the
 convenience of the reader we want to describe the main ideas. Essentially
 this is also the core of Gray's proof of his theorem on Kähler manifolds of
 non-negative sectional curvature.
 
 For any point $p\in M$ let $S(T_pM)$ be the unit sphere in the tangent space
 $T_pM$. Then we have the injective embedding
 $\S^k T^*_pM \rightarrow C^\infty(S(T_pM)), \psi  \mapsto\hat \psi $,
 where the polynomial map $\hat \psi $ is defined by $\hat \psi (X) =
\frac{1}{k!} \,\psi (X, \ldots, X)$. The embedding can be used to define a new scalar
 product  $\tilde g$ on the space of symmetric tensors. For all
 $\psi_1, \psi_2 \in \S^kT^*M$ we define
 $
  \tilde g (\psi_1, \psi_2)(p) : = \int_{S(T_pM)}\hat \psi_1(X)
  \hat \psi_2(X) \mathrm{vol} .
 $
 From Schur's Lemma it follows the existence of a symmetric, positive definite,
 holonomy invariant endomorphism $C$ such that $ \tilde g (\psi_1, \psi_2) =
 g(C \psi_1, \psi_2)$. Note that the eigenspaces of $C$ are just the spaces
 $\S^r_0T^*M$ of trace free symmetric $r$--tensors and that $C$ commutes with
 the curvature endomorphism $q(R)$. Finally a short calculations shows that 
 $
  \tilde g (q(R)\psi_1, \psi_2)  =
  \int_{S(T_pM)} R(X, d\hat\psi_1(X), d\hat\psi_2(X), X) \mathrm{vol}
 $.
 Hence under the assumption of non-negative sectional curvature, it follows
 that $q(R)$ has only non-negative eigenvalues, i.e. $q(R) \ge 0$.
 In the integral above,  $d\hat\psi(X)$ has to be interpreted as the gradient of the function $\hat\psi$
 in the point $X \in S(T_pM)$, i.e. it is again a tangent vector to $M$ in the point $p$.
 \qed
 
 \medskip

 \noindent
 Another fundamental result which we will use below is the following curvature
 formula. It can be proved by a local calculation, but it is also well known
 from the work on the Ricci flow (see e.g. \cite{petersen}, Theorem 9.4.2).

 \begin{Proposition}\label{curv22}
  Let $(\,M,\,g\,)$ be a Riemannian manifold with $\nabla\Ric\,=\,0$.
  Then the Riemannian curvature tensor $R$ considered as a section of
  $\S^2\Lambda^2T^*M$ satisfies:
  $$
   \nabla^*\nabla R
   \;\;=\;\;
   -\;\frac12\,q(R)\,R\ .
  $$
 \end{Proposition}
\subsection{Curvature of Kähler manifolds}
 The holonomy constraint $\nabla J\,=\,0$ imposed on a Kähler manifold
 $(\,M,\,g,\,J\,)$ implies the additional symmetry $[\,R_{X,Y},
 \,J\,]\,=\,0$ of its Riemannian curvature tensor $R$. Using the 
 complex structure $J$ the concept of sectional curvature can be adapted
 to Kähler geometry in form of the holomorphic bisectional curvature
 \begin{equation}\label{hb}
 K(\,X,\,Y\,)\,:=\,R(\,X,\,JX,\,JY,\,Y\,)
 \end{equation} 
 or the holomorphic sectional curvature $S(\,X\,)\,:=\,K(\,X,X\,)$.
 Alternatively the holomorphic sectional curvature can be viewed as
 a symmetric $4$--tensor $S\,\in\,\Gamma( \,\S^4T^*M\,)$ with
 $$
  S(\,X,Y,U,V\,)
  \;\;:=\;\;
  8\;\Big(\;R(X,JY,JU,V\,)
  \;+\;R(\,X,JU,JV,Y\,)
  \;+\;R(\,X,JV,JY,U\,)\;\Big)\ ,
 $$
 these two points of view of $S$ are reconciled by the relation $S(\,X\,)
 \,=\,\frac1{4!}\,S(\,X,X,X,X\,)$. It is possible to reconstruct the
 Riemannian curvature tensor $R$ from $S$, for example via:
 $$
  R(\,X,\,Y,\,U,\,V\,)
  \;\;=\;\;
  \frac1{32}\;\Big(\;S(\,X,\,JY,\,JU,\,V\,)\;-\;S(\,X,\,JY,\,U,\,JV\,)
  \;\Big)\ .
 $$
 In consequence $R\,\longmapsto\,S$ is an injective parallel bundle map
 (see \cite{gregor}, Lemma 2.2). For a Kähler manifold non--negative sectional
 curvature directly implies non--negative holomorphic bisectional curvature
 due to the well--known relation:
 \begin{equation}\label{bihol2}
  K(\;X,\;Y\;)
  \;\;=\;\;
  R(\,X,\,Y,\,Y,\,X\,)\;+\;R(\,X,\,JY,\,JY,\,X\,)\ .
 \end{equation}
 One aim of our article is to show that a similar implication is not true
 for the quaternionic bisectional curvature on quaternion--Kähler manifolds.
\subsection{The new proof of  Gray's theorem}
 Let $(\,M,\,g,\,J\,)$ be a compact Kähler manifold  of non--negative
 sectional curvature and constant scalar curvature. Then Gray's theorem
 follows from two claims, which we prove first.

 \begin{Lemma}
  If the scalar curvature $\kappa$ of a compact Kähler manifold $M$ of
  non--negative sectional curvature is constant, then the Ricci tensor
  is parallel.
 \end{Lemma}

 \proof
 Let $ \S^{2,+} TM$ be the space of symmetric $2$-tensors commuting with $J$.
 Then the map $\theta :  \S^{2,+} TM \rightarrow \Lambda^{1,1} TM$ defined
 by $\theta(h) = h\circ J$
 is a parallel bundle isomorphism. In particular
 we have $\theta(g)\,=\,\omega$, where $\omega$ is the Kähler form
 and $\theta(\Ric) = \rho$, where $\rho$ is the Ricci-form. Recall at
 this point that the Ricci form of a Kähler manifold is closed, namely
 it represents the first Chern class of the holomorphic tangent bundle
 $T^{1,0}M$ up to a constant factor. Since the scalar curvature $\kappa$
 is constant we find $d^*\rho\,=\,0$, as follows from $\delta\Ric\,=\,
 -\frac12\,d\,\kappa$. Hence the Ricci form is harmonic and the universal
 Weitzenböck formula on alternating $2$--forms gives:
 $
  0 = \Delta\rho = \nabla^*\nabla \rho + q(R) \rho\ .
 $
 Composing this with the inverse of the bundle map $\theta$ we obtain the corresponding equation:
 $
  \nabla^*\nabla \Ric = -q(R) \Ric .
 $
 Taking the $L^2$-scalar product with $\Ric$ yields 
 $ 
  (q(R) \Ric, \Ric)_{L^2} = - \| \nabla \Ric \|^2 \le 0 .
 $
 On the other side, since $\Ric$ is a symmetric $2$--tensor and we have
 non--negative sectional curvature by assumption,  Proposition
 \ref{curv11} implies $(q(R) \Ric, \Ric)_{L^2}\,\geq\,0$. Combining the
 two inequalities we conclude that the Ricci tensor has to be parallel. 
 \qed

 \begin{Lemma}
  If the Ricci tensor of a compact Kähler manifold $M$ of non--negative
  sectional curvature is parallel, then the Riemannian curvature tensor
  is parallel as well.
 \end{Lemma}
 
 \proof
 According to Proposition \ref{curv22} the curvature tensor $R$ of a manifold
 $M$ with parallel Ricci tensor satisfies the elliptic partial different
 equation $\nabla^*\nabla R\,=\,-\frac12\,q(R)R$. In turn the image $S$ of
 the curvature tensor $R$ of a Kähler manifold under the parallel bundle
 map $R\longmapsto S$ satisfies the analogous equation $\nabla^*\nabla S
 \,=\,-\frac12\,q(R)S$. Taking the $L^2$--scalar product with $S$ we obtain
 $(\,q(R)S,\,S\,)_{L^2}\,=\,-\,2\,|\!|\,\nabla S\,|\!|^2_{L^2}\,\leq\,0$.
 The assumption of non--negative sectional curvature on the other hand
 tells us $(\,q(R)S,\,S\,)_{L^2}\,\geq\,0$ for the symmetric $4$--tensor
 $S$, and so we conclude that $S$ has to be parallel. Due to injectivity
 of the parallel bundle map $R\longmapsto S$ the curvature tensor $R$ itself has
 to be parallel as well.
 \qed
\section{Curvature of Quaternion--Kähler Manifolds}
\label{three}
 \noindent
 A quaternion--Kähler manifold of quaternionic dimension $n\,\geq\,2$ is a
 Riemannian manifold $(\,M,\,g\,)$ of dimension $m\,=\,4n$ endowed with a
 parallel algebra subbundle $QM\,\subseteq\,\End\,TM$ with respect to the
 Levi--Civita connection $\nabla$. The  fibers of $QM$ are unital subalgebras, i.e.  $\id_{T_pM}\,\in\,Q_pM$, 
 isomorphic to the quaternions
 $Q_pM\,\cong\,\H$ in every point $p\,\in\,M$, such that $g(FX,Y)\,=\,
 g(X,\overline{F}Y)$ for all $X,\,Y\,\in\,T_pM$ and $F\,\in\,Q_pM$. The
 twistor space associated to a quaternion--Kähler manifold $M$ is the
 subbundle of all square roots of $-\id_{T_pM}$ in $QM$:
 $$
  ZM
  \;\;:=\;\;
  \{\;\;(\,p,\,I\,)\;\;|\;\;p\,\in\,M,\;I\,\in\,\,Q_pM\textrm{\ with\ }\,
  I^2\,=\,-\id_{T_pM}\;\;\}
  \;\;\subseteq\;\;
  QM\ .
 $$
 Needless to say the twistor space is a $\C\P^1$--bundle over $M$, more
 precisely $ZM$ equals the unit sphere bundle $ZM\,\subseteq\,\Im\,QM$ of
 the imaginary subbundle $\Im\,QM\,\subseteq\,QM$ with respect to the standard
 scalar product $\gamma(\,F_1,\,F_2\,)\,:=\,\Re(\overline{F}_1F_2)$. According
 to Alekseevskii \cite{alek1} the curvature tensor of a quaternion--Kähler
 manifold of scalar curvature $\kappa$ equals the sum $R\,=\,R^\hyper\,+\,
 R^{\H\P}$ of a curvature tensor $R^\hyper$ of hyperkähler type and the
 curvature tensor $R^{\H\P}$ of the quaternionic projective space of the
 same dimension and scalar curvature, namely
 \begin{equation}\label{riem}
  \begin{array}{lcl}
   R_{X,\,Y}^{\H\P}
   &=&
   -\;{\displaystyle\frac\kappa{16\,n\,(n+2)}}\,
   \Big(\;X\wedge Y\,+\,IX\wedge IY\,+\,JX\wedge JY\,+\,KX\wedge KY
   \\[-2pt]
   &&
   \qquad\qquad\qquad\quad\;\;
   +\;2\,g(IX,Y)\,I\,+\,2\,g(JX,Y)\,J\,+\,2\,g(KX,Y)\,K\;\Big)\ ,
  \end{array}
 \end{equation}
 where $I,\,J,\,K\,\in\,Q_pM$ is a quaternionic triple $I^2\,=\,J^2\,=\,K^2
 \,=\,IJK\,=\,-\id_{T_pM}$ in a given point $p\,\in\,M$ or alternatively an
 oriented orthonormal basis of the imaginary subbundle $\Im\,Q_pM\,\subseteq\,
 Q_pM$ with respect to the standard scalar product $\gamma$.

 The twistor space $ZM$ associated to a quaternion--Kähler manifold $M$ of
 positive scalar curvature $\kappa\,>\,0$ admits an integrable complex
 structure $I^Z$ and a compatible metric $g^Z$ making $ZM$ a Kähler--Einstein
 manifold of positive scalar curvature (\cite{salamon}, Theorem 6.1). In order
 to define $I^Z$ and $g^Z$ we use parallel transport along curves to split the
 tangent space of the twistor space in a point $(\,p,\,I\,)\,\in\,ZM$ into the
 direct sum of $T_pM$ and the subspace $Q^-_pM\,:=\,\span_\R\{\,J,\,K\,\}$ of
 elements anticommuting with $I$ by means of the isomorphism
 \begin{equation}\label{split}
  T_{(\,p,\,I\,)}ZM\;\stackrel\cong\longrightarrow\;T_pM\;\oplus\;Q_p^-M,
  \qquad\left.\frac d{dt}\right|_0(\,p_t,\,I_t\,)\;\longmapsto\;
  \left.\frac d{dt}\right|_0p_t\;\oplus\;\frac12\,I\left(\left.
  \frac\nabla{dt}\right|_0I_t\right)\ ,
 \end{equation}
 where $\left.\frac\nabla{dt}\right|_0I_t$ is the covariant derivative of
 the endomorphism field $I_t\,\in\,Q_{p_t}M$ along the curve $t\longmapsto
 p_t$ in $M$; the sole purpose of the extraneous multiplication with
 $\frac12I$ is to make subsequent formulas more readily comparable with
 the formulas given in \cite{aleksandrov}.

 Under the isomorphism \eqref{split} every tangent vector to the twistor space
 in a given point $(\,p,\,I\,)\,\in\,ZM$ reduces to a sum $X\,\oplus\,x$ with
 $X\,\in\,T_pM$ and $x\,\in\,Q^-_pM$. In this notation the integrable almost
 complex structure $I^Z$ and the compatible metric $g^Z$ can be defined by:
 \begin{eqnarray}
  I^Z(\;X\,\oplus\,x\;)\qquad\quad
  &:=&
  IX\,\oplus\,Ix
  \\[-2pt]
  g^Z(\;X\,\oplus\,x,\;Y\,\oplus\,y\;)
  &:=&
  g(\;X,\;Y\;)\;+\;\frac{16\,n\,(n+2)}\kappa\,\gamma(\;x,\;y\;)\ .
 \end{eqnarray}
 In analogy to the curvature decomposition $R\,=\,R^\hyper\,+\,R^{\H\P}$ for
 $M$ the curvature tensor of the twistor space $ZM$ for the metric $g^Z$
 decomposes into the sum $R^Z\,=\,R^\hyper\,+\,R^{\C\P}$ of the hyperkähler
 part $R^\hyper$ of the curvature of $M$ extended horizontally and the
 curvature tensor $R^{\C\P}$ of the complex projective space of the same
 dimension and scalar curvature:
 \begin{equation}\label{rz}
  R^Z_{X\,\oplus\,x,\,Y\,\oplus\,y}(\,U\,\oplus\,u\,)
  \;\;=\;\;
  (\,R^\hyper_{X,\,Y}U\,\oplus\,0\,)
  \;+\;R^{\C\P}_{X\,\oplus\,x,\,Y\,\oplus\,y}(\,U\,\oplus\,u\,)\ .
 \end{equation}
 Explicitly the complex dimension of the twistor space $ZM$ reads $2n+1$ and
 the equality $\frac\kappa{16n(n+2)}\,=\,\frac{\kappa^Z}{4(2n+1)(2n+2)}$
 relates the scalar curvatures $\kappa$, $\kappa^Z$ of $M$ and $ZM$, in
 turn we find
 \begin{eqnarray}\label{cpr}
   R^{\C\P}_{X\,\oplus\,x,\,Y\,\oplus\,y}
   &=&
   -\;\frac\kappa{16\,n\,(n+2)}\,\Big(\;
   (X \oplus x)\,\wedge\,(Y\oplus y)
   \;+\;I^Z(X \oplus x)\,\wedge\,I^Z(Y\oplus y)
   \\[-10pt]
   \nonumber
   &&
   \hbox to170pt{\hfill}
   +\;2\,g^Z(\,I^Z(X\oplus x),\,Y\oplus y\,)\;I^Z\;\Big) \ .
 \end{eqnarray}
 Intuitively the curvature decomposition \eqref{rz} is a
 direct consequence of the isometry $\C\P^{2n+1}\,\cong\,Z(\H\P^n)$ for the
 quaternionic projective space with respect to the Fubini--Study metrics on
 both $\H\P^n$ and $\C\P^{2n+1}$, a formal demonstration however is
 significantly more involved. Equation \eqref{rz} can be seen as a summary
 of the explicit curvature equations in \cite{aleksandrov} for horizontal
 and vertical tangent vectors.

 An alternatively way to derive equation \eqref{rz} is to reinterprete the
 quaternionic frame bundle of $M$ as a reduction of the frame bundle of the
 twistor space $ZM$. Under this reinterpretation the components of the
 soldering and connection forms for the Levi--Civita connection $\nabla$
 on $M$ reassemble into a soldering and a connection form for $ZM$, which
 endow the twistor space $ZM$ with a metric connection $\overline\nabla^Z$
 with known torsion and curvature tensor. Transgressing from $\overline
 \nabla^Z$ to the Levi--Civita connection $\nabla^Z$ is a rather tedious,
 but straightforward calculation establishing \eqref{rz}. The viability of
 this argument by itself implies that $R^Z$ equals $R$ up to an algebraic
 translation, which is fixed by the special case $M\,=\,\H\P^n$.
 
 \pfill
 Motivated by the work of Chow \& Yang \cite{chow} and the curvature
 decomposition \eqref{rz} we want to study the quaternionic bisectional
 curvature of a quaternion--Kähler manifold $M$ of positive scalar curvature
 $\kappa\,>\,0$ for $I\,\in\,Z_pM$ and tangent vectors $X,\,Y\,\in\,T_pM$ in
 a point $p\,\in\,M$
 $$
  K_{\H,\,I}(\,X,\,Y\,)
  \;\;:=\;\;
  R(\,X,\,IX,\,IY,\,Y\,)\;+\;
  \frac\kappa{8\,n\,(n+2)}\Big(g(\,JX,\,Y\,)^2\;+\;g(\,KX,\,Y\,)^2\Big)\ ,
 $$
 which does not depend on the choice of an orthonormal basis $J,\,K$ for
 $Q^-_pM$.
 
 \begin{Definition}[Quaternionic Bisectional Curvature of Quaternionic Lines]
 \hfill\label{qkb}\break
  The quaternionic bisectional curvature of quaternionic lines $L_1,\,L_2
  \,\subseteq\,T_pM$ in the tangent space of a quaternion--Kähler manifold $M$
  is defined in Chow \& Yang \cite{chow} as the infimum of the quaternionic
  bisectional curvature taken over $I\,\in\,Z_pM$ and unit vectors $X\,\in\,
  L_1,\,Y\,\in\,L_2$:
  $$
   K_\H(\,L_1,\,L_2\,)
   \;\;:=\;\;
   \inf_{\begin{array}{c}\\[-15pt]\scriptstyle I\,\in\,Z_pM,\;X\,\in\,L_1,\;
    Y\,\in\,L_2\\[-3pt]\scriptstyle g(X,X)\,=\,1\,=\,g(Y,Y)
    \end{array}}K_{\H,\,I}(\,X,\,Y\,)\ .
  $$
 \end{Definition}

 \noindent
 The quaternionic bisectional curvature of a quaternion--Kähler manifold $M$
 is closely related to the holomorphic bisectional curvature $K^Z$ of the
 associated twistor space $ZM$:

 \begin{Proposition}[Holomorphic Bisectional Curvature of Twistor Spaces]
 \hfill\label{twistor}\break
  The holomorphic bisectional curvature $K^Z$ of the twistor space $ZM$
  of a quaternion--Kähler manifold $M$ of positive scalar curvature
  $\kappa\,>\,0$ in a point $(\,p,\,I\,)\,\in\,ZM$ equals:
  $$
   K^Z(\,X\oplus x,\,Y\oplus y\,) 
   \;\;=\;\;
   K_{\H,\,I}(\,X,\,Y\,)
   \;+\;K^{\C\P}(\,X\oplus x,\,Y\oplus y\,)
   \;-\;K^{\C\P}(\,X\oplus 0,\,Y\oplus 0\,)\ .
  $$
 \end{Proposition}

 \proof
 On the twistor space the definition \eqref{hb} of the holomorphic bisectional
 curvature reads:
 $$
  K^Z(\,X\oplus x,\,Y\oplus y\,)
  \;\;:=\;\;
  g^Z(\;R^Z_{X\oplus x,\,I^Z(\,X\oplus x\,)}I^Z(\,Y\oplus y\,),\;Y\oplus y\;)
  \ .
 $$
 Due to the curvature decomposition \eqref{rz} and the decomposition $R\,=\,
 R^\hyper\,+\,R^{\H\P}$ of the curvature tensor of the quaternion--Kähler
 manifold $M$ we can expand this definition to read:
 \begin{eqnarray*}
  K^Z(\,X\oplus x,\,Y\oplus y\,)
  &=&
  g(\;R^\hyper_{X,\,IX}IY,\;Y\;)
  \;+\;K^{\C\P}(\;X\oplus x,\;Y\oplus y\;)
  \\
  &=&
  R(\,X,\,IX,\,IY,\,Y\,)
  \;-\;g(\;R^{\H\P}_{X,\,IX}IY,\;Y\;)
  \;+\;K^{\C\P}(\;X\oplus x,\;Y\oplus y\;)\ .
 \end{eqnarray*}
 The description \eqref{cpr} of the curvature tensor $R^{\C\P}$ allows us
 on the other hand to expand the holomorphic sectional curvature of complex
 projective spaces on horizontal tangent vectors:
 \begin{eqnarray*}
  K^{\C\P}(\,X\oplus 0,\,Y\oplus 0\,)
  &=&
  -\;g^Z(\;R^{\C\P}_{X\oplus 0,\,IX\oplus 0}(\,Y\oplus 0\,),\;IY\oplus 0\;)
  \\[2pt]
  &=&
  \hphantom+\;\frac\kappa{8\,n\,(n+2)}\,g
  \Big(\;[\;X\,\wedge\,IX\;+\;g(X,X)\,I\;]\;Y,\;IY\;\Big)
  \\
  &=&
  \hphantom+\;\frac\kappa{8\,n\,(n+2)}\;\Big(\;
  g(X,X)\,g(Y,Y)\;+\;g(X,Y)^2\;+\;g(IX,Y)^2\;\Big)\ .
 \end{eqnarray*}
 In consequence the curvature tensor \eqref{riem} of the quaternionic
 projective space satisfies
 \begin{eqnarray*}
  g(\;R^{\H\P}_{X,\,IX}IY,\;Y\;)
  &=&
  \frac\kappa{8\,n\,(n+2)}\;g\Big(\;
  [\;X\,\wedge\,IX\;-\;JX\,\wedge\,KX\;+\;g(X,X)\,I\;]\,Y,\;IY\;\Big)
  \\
  &=&
  K^{\C\P}(\,X\oplus 0,\,Y\oplus 0\,)\;-\;\frac\kappa{8\,n\,(n+2)}
  \Big(\;g(\,JX,\,Y\,)^2\;+\;g(\,KX,\,Y\,)^2\;\Big)\ .
 \end{eqnarray*}
 for all tangent vectors $X,\,Y\,\in\,T_pM$. Inserting this expression into
 the previous expansion of the holomorphic bisectional curvature
 $K^Z(\,X\oplus x,\,Y\oplus y\,)$ of the twistor space we obtain:
 \begin{eqnarray*}
  K^Z(\,X\oplus x,\,Y\oplus y\,)
  &=&
  K^{\C\P}(\,X\oplus x,\,Y\oplus y\,)
  \;-\;K^{\C\P}(\,X\oplus 0,\,Y\oplus 0\,)
  \\
  &&
  +\;R(\,X,IX,IY,Y\,)
  \;+\;\frac\kappa{8\,n\,(n+2)}\Big(g(JX,Y)^2\,+\,g(KX,Y)^2\Big)
  \qed
 \end{eqnarray*}

 \pfill
 Interestingly the crucial difference $K^{\C\P}(X\oplus x,Y\oplus y)\,-\,
 K^{\C\P}(X\oplus 0,Y\oplus 0)$ appearing in Proposition \ref{twistor}
 turns out to be never negative. In order to verify this assertion we
 use the expansion \eqref{cpr} of the curvature tensor $R^{\C\P}$ of
 complex projective space to obtain:
 \begin{eqnarray*}
  \lefteqn{K^{\C\P}(\,X\oplus x,Y\oplus y\,)}
  \quad
  &&
  \\
  &=&
  -\;g^Z(\;R^{\C\P}_{X\,\oplus\,x,\,I^Z(X\,\oplus\,x)}
  (Y\,\oplus\,y),\;I^Z(Y\,\oplus\,y)\;)
  \\
  &=&
  \frac\kappa{8\,n\,(n+2)}\!\!
  \begin{array}[t]{l}
   \Big[\hphantom{+}\left(\,g(\hphantom IX,Y)
   \,+\,\frac{16n(n+2)}\kappa\gamma(\hphantom Ix,y)\,\right)\!
   \left(\,g(\hphantom IX,Y)
   \,+\,\frac{16n(n+2)}\kappa\gamma(\hphantom Ix,y)\,\right)
   \\[4pt]
   \hphantom{\Big[}+\!
   \left(\,g(IX,Y)\,+\,\frac{16n(n+2)}\kappa\gamma(Ix,y)\,\right)\!
   \left(\,g(IX,Y)\,+\,\frac{16n(n+2)}\kappa\gamma(Ix,y)\,\right)
   \\[4pt]
   \hphantom{\Big[}+\!
   \left(\,g(\,X,\,X\,)\,+\,\frac{16n(n+2)}\kappa\gamma(\,x,x\,)\,\right)\!
   \left(\,g(\,Y,\,Y\,)\,+\,\frac{16n(n+2)}\kappa\gamma(\,y,y\,)\,\right)\Big]
   \ .
  \end{array}
 \end{eqnarray*}
 Cancelling terms in the expansion of the three products on the
 right hand side we arrive at:
 \begin{eqnarray*}
  K^{\C\P}(X\oplus x,Y\oplus y)\,-\,K^{\C\P}(X\oplus 0,Y\oplus 0)
  &=&
  \frac{64\,n\,(n+2)}\kappa\,\gamma(\,x,x\,)\,\gamma(\,y,y\,)
  \\
  &&
  +\;2\,\gamma(x,x)\,g(Y,Y)\;+\;2\,g(X,X)\,\gamma(y,y)
  \\[2pt]
  &&
  +\;4\,\gamma(x,y)\,g(X,Y)\;+\;4\,\gamma(Ix,y)\,g(IX,Y)\ .
 \end{eqnarray*}
 In case $x\,=\,0$ or $X\,=\,0$  the right hand side of this equation
 is a sum of positive terms so that $K^{\C\P}(X\oplus x,Y\oplus y)\,-\,
 K^{\C\P}(X\oplus 0,Y\oplus 0)\,\geq\,0$. Similarly we can expand the norm
 squares $\gamma(y,y)$ and $g(Y,Y)$ under the additional assumption
 $x\,\neq\,0\,\neq\,X$ using Parseval's Equality for the orthogonal
 basis $x,\,Ix$ of $Q_p^-M$ and the orthogonal pair $X,\,IX$
 $$
  \gamma(\,y,\,y\,)
  \;\;=\;\;
  \frac{\gamma(x,y)^2\,+\,\gamma(Ix,y)^2}{\gamma(x,x)}
  \qquad
  g(\,Y,\,Y\,)
  \;\;=\;\;
  g(\,Y^\perp,\,Y^\perp\,)\;+\;
  \frac{g(X,Y)^2\,+\,g(IX,Y)^2}{g(X,X)}
 $$
 with a remainder $Y^\perp\,\in\,\{\,X,\,IX\,\}^\perp$. Hence, 
 the previous difference can be rewritten as 
 \begin{eqnarray*}
  \lefteqn{K^{\C\P}(\,X\oplus x,\,Y\oplus y\,)
   \;-\;K^{\C\P}(\,X\oplus 0,\,Y\oplus 0\,)}
  \quad
  &&
  \\[3pt]
  &=&
  \frac{64\,n\,(n+2)}\kappa\,\gamma(\,x,x\,)\,\gamma(\,y,y\,)
  \;+\;2\,\gamma(\,x,x\,)\,g(\,Y^\perp,Y^\perp\,)
  \\
  &&
  +\;2\,\Big(\,\Lambda^{+\frac12}g(\,X,Y\,)
  \,+\,\Lambda^{-\frac12}\gamma(\,x,y\,)\,\Big)^2
  +\;2\,\Big(\,\Lambda^{+\frac12}g(\,IX,Y\,)
  \,+\,\Lambda^{-\frac12}\gamma(\,Ix,y\,)\,\Big)^2
  \;\;\geq\;\;
  0
 \end{eqnarray*}
 i.e. as a sum of squares with the auxiliary scaling factor $\Lambda\,:=\,
 \frac{\gamma(\,x,\,x\,)}{g(X,X)}$. Combining the non--negativity of the
 difference $K^{\C\P}(X\oplus x,Y\oplus y)\,-\,K^{\C\P}(X\oplus 0,Y\oplus 0)$ 
 for all pairs $X\oplus x,\,Y\oplus y\,\in\,T_{(\,p,\,I\,)}ZM$ of
 tangent vectors with Proposition \ref{twistor} we prove the first part of:
 
 \begin{Corollary}[Non--Negative Quaternionic Bisectional Curvature]
 \hfill\label{nnq}\break
  The quaternionic bisectional curvature of a quaternion--Kähler manifold
  $M$ of positive scalar curvature provides a lower bound for the holomorphic
  bisectional curvature of its associated twistor space $ZM$ in the sense that
  for all tangent vectors $X\oplus x,\,Y\oplus y\,\in\,T_{(\,p,\,I\,)}ZM$:
  $$
   K^Z(\;X\,\oplus\,x,\;Y\,\oplus\,y\;)
   \;\;\geq\;\;
   K_{\H,\,I}(\;X,\,Y\;)\ .
  $$
  In particular, a complete quaternion--Kähler manifold $M$ of positive scalar
  and non--negative quaternionic bisectional curvature is isometric to the
  quaternionic projective space $\H\P^n$.
 \end{Corollary}

 \proof
 Concentrating on the second statement we observe that the twistor space
 $ZM$ associated to a quaternion--Kähler manifold $M$ of positive scalar
 curvature $\kappa\,>\,0$ and non--negative quaternionic bisectional curvature
 $K_\H\,\geq\,0$ has non--negative holomorphic bisectional curvature
 $K^Z\,\geq\,0$ and hence is a locally hermitean symmetric space according
 to a result of Mok \& Zhong \cite{Mok}. In turn the curvature decomposition
 \eqref{rz} of the curvature tensor of $ZM$ implies that the horizontal lift
 of the hyperkähler part $R^\hyper$ of the curvature tensor $R\,=\,R^\hyper
 \,+\,R^{\H\P}$ of $M$ is parallel on $ZM$ and thus a fortiori parallel on $M$.
 In consequence the quaternion--Kähler manifold $M$ is locally symmetric
 and so a Wolf space. Nevertheless the only Wolf space with a locally symmetric
 twistor space is the quaternionic projective space $M\,=\,\H\P^n$ with its
 associated twistor space $ZM\,=\,\C\P^{2n+1}$.
 \qed

\medskip

 \pfill
 Before closing this section we want to discuss the proper analogue for
 quaternion--Kähler manifolds of the important relation \eqref{bihol2}
 between the holomorphic bisectional and the sectional curvature of Kähler
 manifolds. In this proper analogue the quaternionic bisectional curvature
 needs to be replaced by a very similar curvature expression with a different
 scaling factor for the second term, which lacks however a direct geometric
 interpretation in terms of the biholomorphic sectional curvature of
 the associated twistor space $ZM$:

 \begin{Proposition}[Quaternionic Bisectional vs.~Sectional Curvature]
 \hfill\label{qbs}\break
  The quaternionic bisectional curvature of a quaternion--Kähler manifold
  $M$ can not be expressed directly in terms of sectional curvatures of $M$.
  A modified version of the quaternionic bisectional curvature satisfies
  however the relation we expected to see from the Kähler case:
  \begin{eqnarray*}
   \lefteqn{g(\;R_{X,\,IX}IY,\;Y\;)\;+\;\frac\kappa{2\,n\,(n+2)}
   \,\Big(\;g(\,JX,\,Y\,)^2\;+\;g(\,KX,\,Y\,)^2\;\Big)}
   \hbox to180pt{\hfill}
   &&
   \\[-2pt]
   &=&
   g(\,R_{X,\,Y}Y,\,X\,)\;+\;g(\,R_{X,\,IY}IY,\,X\,)\ .
  \end{eqnarray*}
 \end{Proposition}

 \proof
 Recall to begin with that the curvature tensor of a quaternion--Kähler
 manifold $M$ decomposes into the sum $R\,=\,R^\hyper+R^{\H\P}$ of a
 curvature tensor $R^\hyper$ of hyperkähler type and the curvature
 tensor $R^{\H\P}$ of quaternionic projective space. Using the explicit
 description \eqref{riem} of the latter we find for all vectors $X,\,Y\,\in\,
 T_pM$ tangent to $M$
 \begin{eqnarray*}
  [\,R_{X,\,Y},\,I\,]
  \;\;=\;\;
  [\,R^{\H\P}_{X,\,Y},\,I\,]
  &=&
  -\;\frac\kappa{16\,n\,(n+2)}
  \,\Big(\,2\,g(\,JX,\,Y\,)\,[\,J,\,I\,]
  \,+\,2\,g(\,KX,\,Y\,)\,[\,K,\,I\,]\,\Big)
  \\
  &=&
  \hphantom-\;\frac\kappa{4\,n\,(n+2)}
  \,\Big(\;g(\,JX,\,Y\,)\,K\;-\;g(\,KX,\,Y\,)\,J\;\Big)\ ,
 \end{eqnarray*}
 because both $R^\hyper_{X,\,Y}$ and $4\,\pr_{\sp(n)}(\,X\wedge Y\,)
 \,:=\,X\wedge Y+IX\wedge IY+JX\wedge JY+KX\wedge KY$ commute with
 $I,\,J$ and $K$. In light of this commutator formula a standard argument
 leading to relation \eqref{bihol2} for Kähler manifolds can be mended to
 read in the quaternion--Kähler case
 \begin{eqnarray*}
  \lefteqn{g(\;R_{X,\,IX}IY,\;Y\;)}
  &&
  \\[4pt]
  &=&
  g(\;I\,R_{X,\,IY}Y,\;X\;)\;-\;g(\;I\,R_{X,\,Y}\,IY,\;X\;)
  \\
  &=&
  g\Big(\;R_{X,\,IY}IY\,-\,[\,R_{X,\,IY},\,I\,]\,Y,\;X\;\Big)
  \;+\;g\Big(\;R_{X,\,Y}Y\,+\,[\,R_{X,\,Y},\,I\,]\,IY,\;X\;\Big)
  \\
  &=&
  g(\;R_{X,\,IY}IY,\,X\;)\;-\;\frac\kappa{4\,n\,(n+2)}\Big(\,
  g(\,JX,IY\,)\,g(\,KY,X\,)\,-\,g(\,KX,IY\,)\,g(\,JY,X\,)\,\Big)
  \\
  &&
  +\;g(\;R_{X,\,Y}Y,\,X\;)\;+\;\frac\kappa{4\,n\,(n+2)}\Big(\,
  g(\,JX,Y\,)\,g(\,KIY,X\,)\,-\,g(\,KX,Y\,)\,g(\,JIY,X\,)\,\Big)
 \end{eqnarray*}
 where the first equality is just the first Bianchi identity together with the
 skew symmetry of $R_{X,\,IY}$ and $I$. With the due simplifications we arrive
 at the stipulated identity.
 \qed
\section{Wolf spaces}\label{wolf}
 From a close examination of Cartan's list of symmetric spaces, it can be
 deduced that a compact quaternion--Kähler symmetric space exists for each
 type of compact simple Lie group. Wolf's remarkable construction of this
 fact \cite{wolf},  arises from a study of quaternion--Kähler symmetric
 spaces from the  point of view of roots. Wolf devised a method in which
 one begins with a centreless simple Lie group $G,$ its Lie algebra $\g$,
 and a Cartan subalgebra $\t$. By selecting a maximal root, a $\sp(1)$
 subalgebra is singled out. This subalgebra together with its centraliser
 $\k_1$ in $\g$ form the holonomy algebra $\k\,=\,\sp(1)\,\oplus\,\k_1$
 of the quaternion--Kähler symmetric space $M\,=\,G/Sp(1)K_1.$ These
 Wolf spaces exhaust the list of all compact quaternion--Kähler manifolds
 of positive scalar curvature known to date:
 \begin{center}
  \begin{tabular}{|l|l|l|l|}
   \hline
   $G$ & $K$ & $\dim M$  &
   \\
   \hline
   $\SU(n+2)$ & $S(U(n)U(2))$ & $4n\;(n\geq 1)$ & $\SU(2)=\Sp(1)$
   \\
   $\SO(n+4)$ & $\SO(n)\SO(4)$ & $4n\;(n\geq 3)$ & $\SO(4)=\Sp(1)\Sp(1)$
   \\
   $\Sp(n+1)$ & $\Sp(n)\Sp(1)$ & $4n\;(n\geq 1)$ & $\H\P^n$
   \\
   $E_6$ & $\SU(6)\Sp(1)$ & $40$ &
   \\
   $E_7$ & $\Spin(12)\Sp(1)$ & $64$&
   \\
   $E_8$ & $E_7\Sp(1)$ & $112$&
   \\
   $F_4$ & $\Sp(3)\Sp(1)$ & $28$&
   \\
   $G_2$ & $\SO(4)$ & $8$& $\SO(4)=\Sp(1)\Sp(1)$
   \\
   \hline
  \end{tabular}
 \end{center}
 In the non--compact case there exist examples of homogeneous, non--symmetric
 quaternion--Kähler manifolds due to Alekseevskii \cite{alek1,alek2}.
\subsection{Generalities on the compact real form.}
 Let $\g^\C$ be a complex, simple Lie algebra and
 $\h\,\subseteq\,\g^\C$ a Cartan subalgebra with
 associated root system $\Delta\subset\h^*.$ The Killing form
 of $\g^\C$ will be denoted  $B.$ For every root
 $\alpha\,\in\,\Delta$ there exists a unique Lie algebra element
 $H_\alpha\,\in\,\h$ such that $B(H,H_\alpha)\,=\,\alpha(H)$ for all
 $H\,\in\,\h$; this in turn induces a natural scalar product in
 $\Delta$  defined by $\<\alpha,\beta>\,=\,B(H_\alpha,H_\beta)$
 for all $\alpha,\beta\in\Delta$. The root--space decomposition of
 $\g^\C$ is given by
 $$
  \g^\C
  \;\;=\;\;
  \h \oplus\bigoplus_{\alpha\in\Delta}\g^\alpha\ ,
 $$
 where the root spaces $\g^\alpha$ are defined as usual as the common
 eigenspaces of the adjoint action of the Cartan subalgebra
 $\h$ on $\g^\C:$
 $$
  \g^\alpha
  \;\;=\;\;
  \{\;X\in\g^\C\;:\;[H,X]=\alpha(H)X,\;\forall H\in\h\;\}.
 $$
 Let $\g \subset \g^\C$ be a compact real form of $\g^\C$, so by  definition
 $\g^\C\,=\,\g\otimes_\R\C$. The short--hand notation $\g^\C$  will stand for
 the complexification of $\g$. Next we recall the standard construction of a
 compact real form $\g$ of $\g^\C$. The  Cartan subalgebra of $\g$ can be
 defined as
 $$
  \t
  \;\;=\;\;
  \bigoplus_{\alpha\in\Delta}\R\,\left( i H_\alpha\right)
  \;\;\subset\;\;
  \h
 $$ 
 Then, $\h = \t^\C,$ and all roots $\alpha\in\Delta\subset\h^*$ assume
 purely imaginary values on $\t,$ i.e.  $\Delta\subset i\t^*$. Since
 $B(iH_\alpha,iH_\alpha) = -|\alpha|^2$ the  Killing form is negative
 definite on $\t$, thus $\t$ is of compact type. Next we can choose for
 each root $\alpha\in \Delta$ a root vector $X_\alpha\in\g^\alpha$ such that 
 $$ 
  [X_\alpha,X_{-\alpha}] = B(X_\alpha,X_{-\alpha})  H_\alpha =  H_\alpha,
 $$  
 so $B(X_\alpha,X_{-\alpha})=1$. Moreover we define:
 $$
  Y_\alpha
  \;\;=\;\;
  X_\alpha-X_{-\alpha}
  \qquad\qquad
  Z_\alpha
  \;\;=\;\;
  i(X_\alpha+X_{-\alpha})\ .
 $$ 
 Since $B(X_\alpha, X_\beta) =0$ if $\alpha + \beta \neq 0$,  it follows
 that $B(Y_\alpha,Y_\alpha)= B(Z_\alpha,Z_\alpha)=-2, B(Y_\alpha,Z_\alpha)=0$.
 Hence, a compact real form $\g$ of $\g^\C$ is defined by
 $$
  \g
  \;\;=\;\;
  \bigoplus_{\alpha\in\Delta} {\mathbb R}\left( i H_\alpha\right)
  \oplus  \bigoplus_{\alpha\in\Delta} \mathbb R Y_\alpha \oplus
  \R Z_\alpha 
  \;\;=\;\;
  \t \oplus \bigoplus_{\alpha\in\Delta} (\g^\alpha \oplus \g^{-\alpha}\cap\g)
 $$
 Note that $\span_\R\{i H_\alpha, Y_\alpha,  Z_\alpha\}$ is isomorphic
 to the Lie algebra $\sp(1)$. Indeed we have the commutator relations
 $[Y_\alpha,Z_\alpha]=2i H_\alpha,  [i H_\alpha, Y_\alpha]= 2Z_\alpha$
 and $ [i H_\alpha, Z_\alpha]=  - 2 Y_\alpha$. The identification with
 $\sp(1)  = \mathrm {Im} \H =\span\{\,i,\,j,\,k\,\}$ is given by: \;
 $
  i  H_\lambda \longmapsto i,\;Y_\lambda \longmapsto j,\;
  Z_\lambda \longmapsto k.
 $
\subsection{Wolf's theory of the maximal root}
 Under the hypothesis of the previous subsection, a root $\lambda\in\Delta$
 is a maximal root if and only if  $\ad(H_\lambda),$ when restricted to
 $(\g^\lambda \oplus \g^{-\lambda})^\perp$, has only eigenvalues $\pm\frac12
 |\lambda|^2$ and $0$. Equivalently, if and only if 
 $$
  \frac{2\<\alpha,\lambda>}{\<\lambda,\lambda>}
  \;\;=\;\;
  \{\,-1,\,0,\,+1\,\}
  \qquad
  \textrm{for all} \;\; \alpha\neq \pm\lambda.
 $$
 Indeed if $\lambda$ is maximal for some ordering and $\alpha\in\Delta^+$
 is a positive root different from $\lambda$, then by studying the $\lambda$--string through
 $\alpha,$ neither $\alpha+\lambda$ nor $\alpha-2\lambda$ are roots,
 but $\alpha-\lambda$ might be a root, thus the scalar product
 $\frac{2\,\<\alpha,\lambda>}{\<\lambda,\lambda>}\,=\,\{0,\pm 1\}.$ 

 From now on we will call such a root $\lambda$ a Wolf root. By virtue
 of the existence and properties of a Wolf root $\lambda$, the complex
 Lie algebra $\g^\C$ has a parabolic splitting
 $\g^\C\,=\,\k^\C \oplus \m^\C$ such that
 $$
  \g^\C
  \;\;=\;\; 
  \underbrace{\left[ (\t^\C
  \oplus\bigoplus_{\alpha\in\lambda^\perp}\g^\alpha )
  \oplus \g^\lambda \oplus\g^{-\lambda}\right]}_{\k^\C}
  \oplus
  \underbrace{\bigoplus_{\<\alpha,\lambda>=\frac12 |\lambda|^2} 
  \g^\alpha\oplus \g^{-\alpha}}_{\m^\C}.
 $$
 Let $\k\,=\,\k^\C\cap\g$ and $\m\,=\,\m^\C \cap \g$. Then $\g = \k \oplus \m$
 and we have the inclusions $[\m, \m] \subset \k$ and $[\k, \m] \subset \m$.
 Hence $(\g, \k) $ is a symmetric pair with corresponding symmetric space
 $M\,=\,G/K$ of compact type and with isotropy representation $\m$. As we
 will recall below, a Wolf root $\lambda$ defines on $M$ the structure of a
 quaternion--Kähler manifold, the so--called Wolf space.
\subsection{The complex structure}
 It is possible to define a complex structure on $\m^\C \subset \g^\C$
 as follows: Let $I\in \t$ be defined as
 $$
  I
  \;\;=\;\;
  i\,\frac{2 H_\lambda}{\<\lambda,\lambda>}
 $$ 
 Hence $\lambda(I)=2i,\,\alpha(I)=0,$ for all $\alpha\in\lambda^\perp.$ We
 consider $I$ as a endomorphism of $\g^\C$ via the adjoint action and write
 $IX=[I,X].$ Then $I$ has eigenvalues $\pm i$ on $\m^\C$ and we can define
 $$
  \m^{1,0}
  \;\;=\;\;
  \{X\in\g^\C \;:\; [I,X]=iX\}
  \;\;=\;\;
  \bigoplus_{\<\alpha,\lambda>=\frac12|\lambda|^2}\g^\alpha,
 $$ 
 $$
  \m^{0,1}
  \;\;=\;\;
  \{X\in\g^\C \;:\; [I,X]=-iX\}
  \;\;=\;\;
  \bigoplus_{\<\alpha,\lambda>=\frac12| \lambda|^2}\g^{-\alpha}\ ,
 $$ 
 so that $\m^{1,0}\oplus \m^{0,1}=\m^\C$. Indeed, for any root $\alpha$ the
 condition  $\<\alpha,\lambda>\,=\,\frac12\<\lambda,\lambda>$ holds, if and
 only if  $\alpha$ can be written as $\alpha = \frac12 \lambda + \alpha_0$,
 where $\alpha_0 $ is orthogonal to $\lambda$. By definition of $I$ this is
 equivalent to $\alpha(I) = i$ and thus to  $[I, X_\alpha] = i X_\alpha$.

 The  action of $I $ on $\m = (\m^{1,0} \oplus \m^{0,1}) \cap \g $ can be
 given  explicitly on the basis vectors $Y_\alpha$ and $Z_\alpha$. For any
 root $\alpha$ with $\<\alpha,\lambda>\,=\,\frac12|\,\lambda\,|^2$ we have
 $$
  \begin{array}{lclclcl} 
   IY_\alpha
   &=&
   [\;I,\;Y_\alpha\;]
   &=&
   (+i)\,X_\alpha\;  - \;(-i)\,X_{-\alpha}
   &=&
   \hphantom{+}\;Z_\alpha,
   \\
   IZ_\alpha
   &=&
   [\;I,\;Z_\alpha\;]
   &=&
   (+i)\,i\,X_\alpha\;  +  \;(-i)\,i\,X_{-\alpha}
   &=&
   -\;Y_\alpha,
  \end{array}
 $$
 thus implying the fundamental relation $I^2\,=\,-\id_\m$ of a complex
 structure on $\m$.
\subsection{The quaternionic structure.}
 The Lie algebra $\sp(1) \cong \mathrm{ span}_{\mathbb R}\{i H_\lambda,
 Y_\lambda, Z_\lambda\}$ acts via the adjoint representation on $\m$.
 More precisely, denoting $ L_\alpha =( \g^\alpha\oplus\g^{-\alpha})\cap \g$
 for  $\alpha \in \Delta$ with $2\<\alpha,\lambda> = |\lambda|^2$ and
 $J ={\rm ad}(Y_\lambda),\; K\,=\,\ad(Z_\lambda)$, we find the commutator
 relations
 $$
  I(L_\alpha) \subset L_\alpha
  \qquad\qquad
  J(L_\alpha) \subset L_{\lambda - \alpha}
  \qquad\qquad
  K(L_\alpha) \subset L_{\lambda -\alpha}\ .
 $$
 The first inclusion follows from $I(X_\alpha) = i X_\alpha$. For the second
 inclusion we note that for any positive root $\alpha$ we have $J(X_\alpha)
 \,=\,[X_\lambda - X_{-\lambda}, X_\alpha]\,\in\,\g^{ - \lambda + \alpha}$,
 because $\alpha + \lambda$ is not a root, as $\lambda$ is maximal and
 similarly we conclude $J(X_{-\alpha}) \in \g^{\lambda  - \alpha}$. The
 last inclusion follows in the same way. In particular we have a non trival
 $\sp(1)$ representation on the real $4$-dimensional space $L_\alpha \oplus
 L_{\lambda-\alpha}$ which has to be the standard action on $\H$. Note that
 $\lambda - \alpha \neq \alpha$ and again $2\<\lambda\,-\,\alpha,\lambda>
 \,=\,|\lambda|^2$. Up to rotation we thus obtain $J ^2\,=\,K ^2\,=\,-\id$
 and $JK\,=\,I$ on the isotropy representation $\m$. This defines a
 quaternionic structure $QM$ on the Wolf space $M\,=\,G/K$, 
 as introduced in the beginning of Section \ref{three}.
 
 We will call $L_\alpha \oplus L_{\lambda-\alpha} \subset \m = T_oM$
 the \textit{quaternionic line} associated to $\alpha\in\Delta$. It is
 the $4$--dimensional subspace of $T_oM$ spanned by $X, IX, JX, KX$
 for any tangent vector $X \in L_\alpha$.
\subsection{Quaternionic bisectional curvature of Wolf spaces}
 Let $M\,=\,G/K$ be a Wolf space of compact type with Cartan decomposition
 $\g\,=\,\k\oplus\m$. On $M$ we fix a Riemannian metric $g$ induced by a
 negative multiple of the Killing form $B$ of $\g$ restricted to the isotropy
 representation $\m$. As a symmetric space the Riemannian curvature of
 $(M, g)$ is given by $g(\,R_{X,\,Y}U,\,V\,)\,=\,-g(\,[X,Y],\,[U,V]\,)$
 for any tangent vectors $X,\,Y,\,U,\,V\,\in\,\m\,=\,T_oM$.

 \begin{Proposition}\label{curv}
  Let $(\,M,\,g\,)$ be a Wolf space of compact type with Wolf root $\lambda$.
  Then for any tangent vectors $X \in L_\alpha$ and $Y \in L_\beta$ we have
  $$
   g(R_{X,IX}IY,Y)
   \;\;=\;\;
   \frac\kappa{2n(n+2)}\frac{\<\alpha,\beta>}{\<\lambda,\lambda>}
   \,g(X,X)\,g(Y,Y)\ .
  $$
 \end{Proposition}

 \proof
 Any vector $X \in L_\alpha$ we can write as $X = X^{1,0} + X^{0,1}$ with
 $X^{1,0} \in \g^\alpha \subset \m^{1,0}$ and $X^{0,1} \in \g^{-\alpha}
 \subset \m^{0,1}$. In particular, we have $g(X,X) = 2g(X^{1,0}, X^{0,1})$,
 where we use the $\C$-bilinear extension to $\m^\C$ of the metric $g$. Then
 $$
  [X, IX]
  \;\;=\;\;
  [X, [I, X]]
  \;\;=\;\;
  [X^{1,0} + X^{0,1}, i  X^{1,0} - i  X^{0,1}]
  \;\;=\;\;
  -2i [X^{1,0},X^{0,1}] \in \t
 $$
 Taking the scalar product with an arbitrary element $H\,\in\,\t$ we find
 \begin{eqnarray*}
  g(i[X^{1,0},X^{0,1}], H)
  &=&
  - ig( X^{1,0}, [H, X^{0,1}])
  \;\;=\;\;
  i \alpha(H) g(X^{1,0},X^{0,1})
  \;\;=\;\;
  \tfrac12  i \alpha(H) g(X, X)
  \\[.5ex]
  &=&
  \tfrac12  g(X, X) g( i \alpha^\sharp, H ) ,
 \end{eqnarray*}
 where we define $\alpha^\sharp$ by $ i \alpha(H) =  g( i \alpha^\sharp, H )$
 for all $H\,\in\,\t$. Note that $\alpha^\sharp\,\in\,i\t$ because of
 $\alpha\,\in\,i\t^*$. It follows
 $
  i\,[\;X^{1,0},\;X^{0,1}\;]
  \,=\,
  \frac12 g(X,X) i \alpha^\sharp
 $.
 Substituting this into the formula for the Riemannian curvature of the
 symmetric metric $g$ we obtain
 \begin{eqnarray*}
  g(R_{X, IX} IY, Y)
  &=&
  g([X, IX], [Y, IY])
  \;\;=\;\;
  4 g(i [X^{1,0},X^{0,1}],i [Y^{1,0},Y^{0,1}])
  \\[1ex]
  &=&
  - g(X,X) g(Y, Y) g(\alpha^\sharp, \beta^\sharp)
  \;\;=\;\;
  -\,g(X,X)\,g(Y, Y)\,\<\alpha, \beta>_c
 \end{eqnarray*}
 where $\<\cdot, \cdot>_c$ is a constant multiple of the scalar product
 $\<\cdot,\cdot>$ induced by the Killing form on the root space $\Delta$.
 In order to remove the dependence of a scaling factor we recall a formula
 for the scalar curvature $\kappa$ in terms of the  length of the Wolf root
 $\lambda$. In (\cite{semmelmann}, Eq. (4.1)) it is proved that:
 $
  \kappa = - 2n(n+2 )\<\lambda,\lambda>_c\ . 
 $
 Hence, we have
 \begin{eqnarray*}
  \tfrac\kappa{2n(n+2)}\tfrac{\<\alpha,\beta>}{\<\lambda,\lambda>}
  \,g(X,X)\,g(Y,Y)
  &=&
  \tfrac\kappa{2n(n+2)}\tfrac{\<\alpha,\beta>_c}{\<\lambda,\lambda>_c}
  \,g(X,X)\,g(Y,Y)
  \\
  &=&
  -\;\<\alpha,\beta>_c\,g(X,X)\,g(Y, Y)\ .
 \end{eqnarray*}
 This then proves the stated curvature formula.
 \qed

 \begin{Corollary}[Quaternionic Bisectional Curvature of Wolf Spaces]
 \hfill\label{qkcurv}\break
  Consider the Wolf space associated to a Wolf root $\lambda$ of a simple
  Lie algebra $\g\,\not\cong\,\sp(1)$. For every tangent vectors $X\,\in\,
  L_\alpha$ and $Y\,\in\,L_\beta$ the quaternionic bisectional curvature
  equals:
  $$
   K_{\H,\,I}(\;X,\;Y\;)
   \;\;=\;\;
   \frac\kappa{2\,n\,(n+2)}\,\left(\;\frac{\<\alpha,\beta>}{\<\lambda,\lambda>}
   \;+\;\frac{\delta_{\beta,\,\lambda-\alpha}}{4}\;\right)
   \;g(\,X,\,X\,)\;g(\,Y,\,Y\,)\ .
  $$
 \end{Corollary}
 
 \proof
 We have to determine the additional summand in the quaternionic bisectional
 curvature, i.e.~we have to compute $g(JX, Y)^2 + g(KX, Y)^2$. Note that
 $L_\alpha \perp L_\beta$ for $\beta \neq \pm \alpha$ and that for any vector
 $X \in L_\alpha $ the space $L_{\lambda - \alpha}$ is spanned by $JX$ and
 $KX$. Hence, if $\pi$ denotes the orthogonal projection of a vector $Y$ to 
 $L_{\lambda - \alpha}$, then $\pi(Y) |X|^2 = g(JX,Y)JX + g(KX,Y)KX$ and
 $\pi(Y)=\delta_{\beta,\lambda-\alpha}Y$. After taking the scalar product
 with $Y$ it follows that
 $$
  \delta_{\beta,\lambda-\alpha} \, |Y|^2 |X|^2
  \;\;=\;\;
  g(\pi(Y), Y) |X|^2
  \;\;=\;\;
  g(JX, Y)^2 + g(KX, Y)^2  
 $$
 From this expression and  Proposition \ref{curv} the formula for
 $K_\H (X, Y)$ immediately follows.
 \qed

 \begin{Theorem}[Negative Quaternionic Bisectional Curvature]
 \hfill\label{main}\break
  For every Wolf space $M$ of compact type different from the quaternionic
  projective space $\H\P^n$ the quaternionic bisectional curvature $K_\H$
  assumes negative values. More precisely for every Wolf root $\lambda$ of
  a simple compact Lie algebra $\g$ different from $\sp(n+1),\,n\,\geq\,0,$
  there exists a long root $\alpha$ satisfying $|\,\lambda\,|^2\,=\,
  2\<\alpha,\lambda>\,=\,|\,\alpha\,|^2$. In consequence the quaternionic
  bisectional curvature becomes strictly negative on all non--zero tangent
  vectors $X\,\in\,L_\alpha$ and $Y\,\in\,L_{\lambda-\alpha}$:
  $$
   K_{\H,\,I}(\;X,\;Y\;)
   \;\;=\;\;
   -\;\frac\kappa{8\,n\,(n+2)}\;g(\,X,\,X\,)\;g(\,Y,\,Y\,)
   \;\;<\;\;
   0\ .
  $$
 \end{Theorem}
 \noindent
 The constraints $|\,\lambda\,|^2\,=\,2\<\alpha,\lambda>\,=\,|\,\alpha\,|^2$
 characterizing the long root $\alpha$ imply of course:
 $$
  \frac{\<\alpha,\lambda-\alpha>}{\<\lambda,\lambda>}
  \;+\;\frac{\delta_{\lambda-\alpha,\,\lambda-\alpha}}4 
  \;\;=\;\;
  \frac12\;-\;1\;+\;\frac14
  \;\;=\;\;
  -\;\frac14\ .
 $$
 Hence the quaternionic bisectional curvature $K_{\H,\,I}(\,X,\,Y\,)\,=\,-
 \frac\kappa{8n(n+2)}\,g(X,X)\,g(Y,Y)$ is strictly negative for non--zero
 tangent vectors $X\,\in\,L_\alpha$, $Y\,\in\,L_{\lambda-\alpha}$ by
 Corollary \ref{qkcurv}. In turn the proof of Theorem \ref{main} reduces to
 the existence of such a long root $\alpha$, or equivalently to:

 \begin{Lemma}[Existence of Pairs of Non--Orthogonal Long Roots]
 \hfill\label{nort}\break
  Every irreducible root system $\Delta$ of type different to $A_1$ and $C_n$
  for some $n\,\geq\,2$ contains a pair $\lambda,\,\alpha\,\in\,\Delta$ of
  linearly independent,  i.e.  $\lambda\,\neq\,\pm\alpha$, 
  non--orthogonal
 long roots.
 \end{Lemma}

 \proof
 Recall to begin with that an irreducible root system $\Delta$ is called
 simply laced provided all its roots have the same length and are thus long
 roots by definition. Our case by case proof begins with a simply laced root
 system $\Delta$ of type $A_n$ and rank $n\,\geq\,1$, which reads
 $$
  \Delta
  \;\;=\;\;
  \{\;\;\pm\,(\,\e_\mu\,-\,\e_\nu\,)\;\;|
  \;\;0\,\leq\,\mu\,<\,\nu\,\leq\,n\;\;\}
 $$
 in terms of weights $\e_0,\,\ldots,\,\e_n$ subject to the linear
 dependency $\e_0+\ldots+\e_n\,=\,0$ and
 $$
  \<\;\e_\mu,\;\e_\nu\;>
  \;\;:=\;\;
  \delta_{\mu\nu}\;-\;\frac1{n+1}
 $$
 up to a common proportionality constant. Excluding the root system $A_1$ by
 requiring $n\,\geq\,2$ we see that the roots $\lambda\,:=\,\e_0-\e_1$ and
 $\alpha\,:=\,\e_0-\e_2$ are not orthogonal due to $\<\lambda,\alpha>\,=\,1$.
 In turn the long roots in a root system $\Delta$ of type $B_n$ or $D_n$ and
 rank $n\,\geq\,3$ can be written
 $$
  \Delta
  \;\;\supseteq\;\;
  \{\;\;\pm\,\e_\mu\,\pm\,\e_\nu\;\;|\;\;1\,\leq\,\mu\,<\,\nu\,\leq\,n\;\;\}
 $$ 
 in terms of an orthogonal basis $\e_1,\,\ldots,\,\e_n$ of weights with
 $\<\e_\mu,\e_\nu>\,=\,\delta_{\mu\nu}$ up to a common proportionality constant.
 All these root systems $\Delta$ contain the pair $\lambda\,:=\,\e_1+\e_2$ and
 $\alpha\,:=\,\e_1+\e_3$ of non--orthogonal long roots with $\<\lambda,\alpha>
 \,=\,1$. In passing we observe that the short roots in root systems of type
 $B_n$ read $\pm\,\e_1,\,\ldots,\,\pm\,\e_n$, while root systems of type $D_n$
 are simply laced. The additional rank assumption $n\,\geq\,3$ in this case is
 justified by the exceptional Dynkin diagram isomorphisms $B_2\,\cong\,C_2$
 and $D_2\,\cong\,A_1\times A_1$.

 Coming to the irreducible root systems of exceptional type we observe that
 irreducible root systems $\Delta$ of type $G_2$ and $F_4$ contain simply
 laced irreducible subsystems of type $A_2$ and $D_4$ respectively consisting
 entirely of long roots, while the simply laced irreducible root systems
 $\Delta$ of exceptional type $E_6,\,E_7$ and $E_8$ contain irreducible
 subsystems of type $A_5,\,A_6$ and $A_7$. Arguing last but not least to
 the contrary we recall that the long roots in an irreducible root system
 $\Delta$ of type $C_n$ and rank $n\,\geq\,2$ read $\pm\,2\e_1,\,\ldots,
 \,\pm\,2\e_n$ for an orthogonal basis $\e_1,\,\ldots,\,\e_n$ of weights
 with scalar product proportional to $\<\e_\mu,\e_\nu>\,=\,\delta_{\mu\nu}$.
 Every linearly independent pair of long roots in $\Delta$ is thus
 automatically orthogonal.
 \qed
\end{document}